\date{November 15, 2006}

\documentclass[11pt,leqno]{article}
\usepackage{graphicx}
\usepackage{amssymb}
\newtheorem{definition}{Definition}[section]
\newtheorem{lemma}[definition]{Lemma}
\newtheorem{theorem}[definition]{Theorem}
\newtheorem{proposition}[definition]{Proposition}

\newtheorem{remark}[definition]{Remark}
\newenvironment{proof}{\noindent{\bf Proof }}{\hfill $\Box$\medskip}

\renewcommand{\a}{\alpha}

\newcommand{\PP}{{\mathbb{P}}}
\newcommand{\ZZ}{{\mathbb{Z}}}
\newcommand{\CC}{{\mathbb{C}}}

\newcommand{\RR}{{\mathbb{R}}}

\newcommand{\cO}{{\mathcal{O}}}

\newcommand{\iso}{\cong}

\newcommand{\Hom}{\hbox{Hom}}
\newcommand{\GL}{\hbox{GL}}
\renewcommand{\Im}{\hbox{Im}\,}
\renewcommand{\Re}{\hbox{Re}}

\newcommand{\surj}{\twoheadrightarrow}
\newcommand{\inc}{\hookrightarrow}
\newcommand{\too}{\longrightarrow}

\newcommand{\la}{\langle}
\newcommand{\ra}{\rangle}

\newcommand{\bV}{\bigwedge V}
\newcommand{\bW}{\bigwedge W}

\begin{document}

\title{Symplectic nilmanifolds with a symplectic non-free $\ZZ_3$-action}

\author{Marisa Fern\'andez and Vicente Mu\~noz}

\maketitle

\begin{abstract}
In \cite{FM4} the authors have introduced a new technique to
produce symplectic manifolds. It consists on taking a symplectic
non-free action of a finite group on a symplectic manifold and
resolving symplectically the singularities of the quotient. This
has allowed to produce the first example of a non-formal simply
connected compact symplectic manifold of dimension $8$. Here we
present another description of such a manifold and we expand on
some of the details concerning its properties.
\end{abstract}

\section{Introduction} \label{sec:introduction}

In \cite{FM4}, the authors have produced the first example of a
simply connected compact symplectic manifold of dimension $8$
which is non-formal.

In general, simply connected compact manifolds of dimension less
than or equal to $6$ are formal \cite{NM,FM1}, and there are
simply connected compact manifolds of dimension greater than or
equal to $7$ which are non-formal
\cite{Oprea,FM2,Rudyak,Cav1,FM3}. This is a problem that can be
tackled by using minimal models \cite{DGMS} and suitable
constructions of differentiable manifolds.

However, if we consider symplectic manifolds, the story is not so
straightforward, basically due to the fact that there are not so
many constructions of symplectic manifolds. In \cite{BT1,BT2}
Babenko and Taimanov give examples of non-formal simply connected
compact {\em symplectic} manifolds of any dimension bigger than or
equal to $10$, by using the symplectic blow-up \cite{McDuff}. They
raise the question of the existence of non-formal simply connected
compact symplectic manifolds of dimension $8$. Examples of these
cannot be constructed by means of symplectic blow-ups. Other
methods of construction of symplectic manifolds, like the
connected sum along codimension two submanifolds \cite{Go}, or
symplectic fibrations \cite {McDuff-Salamon,Th,We} have not
produced such examples so far.

The solution to this question presented in \cite{FM4} uses a new
and simple method of construction of symplectic manifolds. This
method consists on taking quotients of symplectic manifolds by
finite groups and resolving symplectically the singularities.
Starting with a suitable compact non-formal nilmanifold of
dimension $8$, on which the finite group $\ZZ_3$ acts with simply
connected quotient, one gets a simply connected compact symplectic
non-formal $8$--manifold.

In this note, we expand on some of the issues touched in
\cite{FM4}. First we present an alternative description of the
manifold constructed in \cite{FM4}, by using real Lie groups
instead of complex Lie groups (see Section
\ref{sec:preliminaries}). Actually this is the way in which we
first obtained the example; introducing complex Lie groups was an
ulterior simplification. The reason for our choice of symplectic
$8$--dimensional nilmanifold $M$ becomes transparent with the
description that we give in Section \ref{sec:hatE}: it is the
simplest case in which the group $\ZZ_3$ acts not having any
invariant part in the cohomology of degree $1$. In this way we
have a chance to get a simply connected symplectic orbifold
$\widehat M=M/\ZZ_3$, as we prove later it is the case with our
particular choice of $M$ and $\ZZ_3$-action. To get a {\em
smooth\/}  $8$--dimensional symplectic manifold, we have to
resolve symplectically the singularities. For this, in Section
\ref{sec:mainresult}, we take suitable K\"ahler models around each
singular point. It is clear that this method can be used in much
greater generality.

The last issue concerns with the non-formality of the constructed
manifold. Our example of symplectic $8$--manifold has vanishing
odd Betti numbers, therefore its triple Massey products are zero.
Thus the natural way to prove non-formality is to produce the
minimal model of $\widehat M$, but this can be a lengthy task for
large Betti numbers.

In \cite{FM1} the concept of formality is extended to a weaker
notion named as $s$--formality. We shall not review this notion
here, but we want to mention that it has always been a guidance
for us when trying to write down an obstruction to detect
non-formality. Actually, when we spelt out the condition for
$3$--formality, we realised that there is an easily described new
type of obstruction to formality constructed with differential
forms, in spirit similar to Massey products, and which we
christen here
as {\em G-Massey product} (see Definition \ref{def:GMassey} in Section
\ref{sec:formality}). Then, in Section \ref{sec:non-formal}, the
non-formality of $\widehat M$ is easily checked via a non-trivial
G-Massey product. There, we also see that a suitable quadruple
Massey product of $\widehat M$ is non-trivial, although the proof
is definitely more obscure. So we have decided to include both
proofs of the non-formality of $\widehat M$. It would be
interesting to find a space with non-trivial G-Massey products but
with all multiple Massey products trivial.

\bigskip

\noindent {\bf Acknowledgments.} We would like to thank Luis C. de
Andr\'es, Dominic Joyce, Gil Cavalcanti, Ignasi Mundet, John Oprea
and Daniel Tanr\'e for conversations. This work has been partially
supported through grants MCyT (Spain) MTM2004-07090-C03-01,
MTM2005-08757-C04-02  and Project UPV 00127.310-E-15909/2004.

\medskip
\section{Formality and G-Massey products} \label{sec:formality}

Before we start with the construction of symplectic manifolds, we
shall briefly review the notion of formality \cite{FM1, DGMS},
and we shall introduce the {\em G-Massey products}
as an obstruction to this property.
Let $(A,d)$ be a {\it differential graded commutative algebra} over the field
$\RR$ of real numbers. Then $(A,d)$ is said to be {\it minimal\/}
if:
\begin{enumerate}
 \item $A$ is free as an algebra, that is, $A$ is the free
 algebra $\bigwedge V$ over a graded vector space $V=\oplus V^i$, and
 \item there exists a collection of generators $\{ a_\tau,
 \tau\in I\}$, for some well ordered index set $I$, such that
 $\deg(a_\mu)\leq \deg(a_\tau)$ if $\mu < \tau$ and each $d
 a_\tau$ is expressed in terms of preceding $a_\mu$ ($\mu<\tau$).
 This implies that $da_\tau$ does not have a linear part, i.e., it
 lives in $\bV^{>0} \cdot \bV^{>0} \subset \bV$.
\end{enumerate}

Given a differential algebra $(A,d)$, we denote by $H^*(A)$ its
cohomology. $(A,d)$ is connected if $H^0(A)=\RR$. We shall say
that $({\cal M},d)$ is a {\it minimal model} of the differential
algebra $(A,d)$ if $({\cal M},d)$ is minimal and there exists a
morphism of differential graded algebras $\rho\colon {({\cal
M},d)}\longrightarrow {(A,d)}$ inducing an isomorphism
$\rho^*\colon H^*({\cal M})\longrightarrow H^*(A)$ on cohomology.
In~\cite{H} Halperin proved that any connected differential
algebra $(A,d)$ has a minimal model unique up to isomorphism.

A {\it minimal model\/} of a connected differentiable manifold $X$
is a minimal model $(\bigwedge V,d)$ for the de Rham complex
$(\Omega X,d)$ of differential forms on $X$. If $X$ is a simply
connected manifold, then the dual of the real homotopy vector
space $\pi_i(X)\otimes \RR$ is isomorphic to $V^i$ for any $i$.
This relation also happens when $i>1$ and $X$ is nilpotent, that
is, the fundamental group $\pi_1(X)$ is nilpotent and its action
on $\pi_j(X)$ is nilpotent for $j>1$ (see~\cite{DGMS}).

A minimal model $({\cal M},d)$ is said to be {\it formal} if there
is a morphism of differential algebras $\psi\colon {({\cal
M},d)}\longrightarrow (H^*({\cal M}),d=0)$ that induces the
identity on cohomology. We shall say that $X$ is {\it formal\/} if
its minimal model is formal or, equivalently, the differential
algebras $(\Omega X,d)$ and $(H^*(X),d=0)$ have the same minimal
model. Therefore, if $X$ is formal and simply connected, then the
real homotopy groups $\pi_i(X)\otimes \RR$ are obtained from the
minimal model of $(H^*(X),d=0)$.

Many examples of formal manifolds are known: spheres, projective
spaces, compact Lie groups, homogeneous spaces, flag manifolds,
and compact K\"ahler manifolds. The importance of formality in
symplectic geometry stems from the fact that it allows to
distinguish between symplectic manifolds which admit K\"ahler
structures and some which do not \cite{TO}.

In order to detect non-formality, instead of computing the minimal
model, which usually is a lengthy process, we can use Massey
products, which are obstructions to formality. Let us recall its
definition. The simplest type of Massey product is the triple
(also known as ordinary) Massey product. Let $X$ be a (not
necessarily simply connected) manifold and let $a_i \in
H^{p_i}(X)$, $1 \leq i\leq 3$, be three cohomology classes such
that $a_1\cup a_2=0$ and $a_2\cup a_3=0$. The (triple) Massey
product of the classes $a_i$ is defined as the set
  $$
  \langle a_1,a_2,a_3 \rangle  = \{
  [ \alpha_1 \wedge \eta+(-1)^{p_1+1} \xi
  \wedge \alpha_3] \ | \ a_i=[\alpha_i],\ \alpha_1\wedge \alpha_2= d \xi,
  \ \alpha_2\wedge \alpha_3=d \eta \} \subset
  H^{p_1+p_2+p_3-1}(X) \ .
  $$
(The same set is obtained if we fix the forms $\alpha_i$ such that
$a_i=[\alpha_i]$ and we only let $\eta$ and $\xi$ vary.) It is
easily seen that $\langle a_1,a_2,a_3 \rangle$ is a set of the
form $b+ (a_1 \cup H^{p_2+p_3-1}(X) + H^{p_1+p_2-1}(X)\cup a_3)$,
so it gives a well-defined element in
  $$
  \frac{H^{p_1+p_2+p_3-1}(X)}{a_1
  \cup H^{p_2+p_3-1}(X) + H^{p_1+p_2-1}(X)\cup a_3} \; .
  $$
We say that $\langle a_1,a_2,a_3 \rangle$ is trivial if $0\in
\langle a_1,a_2,a_3 \rangle$.

The definition of higher Massey products is as follows
(see~\cite{K,Mas,Ta}). The Massey product $\la
a_1,a_2,\dots,a_t\ra$, $a_i\in H^{p_i}(X)$, $1\leq i\leq t$,
$t\geq 3$, is defined if there are differential forms $\a_{i,j}$
on $X$, with $1\leq i\leq j\leq t$, except for the case
$(i,j)=(1,t)$, such that
 \begin{equation}\label{eqn:gm}
 a_i=[\a_{i,i}], \qquad
 d\,\a_{i,j}= \sum\limits_{k=i}^{j-1} {\bar \a}_{i,k}\wedge
 \a_{k+1,j},
 \end{equation}
where $\bar \a=(-1)^{\deg(\a)} \a$. Then the Massey product is
 $$
 \la a_1,a_2,\dots,a_t \ra =\left\{
 \left[\sum\limits_{k=1}^{t-1} {\bar \a}_{1,k} \wedge
 \a_{k+1,t}\right] \ | \ \a_{i,j} \hbox{ as in (\ref{eqn:gm})}\right\}
 \subset H^{p_1+ \cdots +p_t
 -(t-2)}(X)\, .
 $$
We say that the Massey product is trivial if $0\in \la
a_1,a_2,\dots,a_t\ra$. Note that for $\la a_1,a_2,\dots,a_t\ra$ to
be defined it is necessary that $\la a_1,\dots,a_{t-1}\ra$ and
$\la a_2,\dots,a_t\ra$ are defined and trivial.

The existence of a non-trivial Massey product is an obstruction to
formality. Concretely, we have the following result, for whose
proof we refer to~\cite{DGMS,Ta}.

\begin{lemma} \label{lem:criterio1}
 If $X$ has a non-trivial Massey product then $X$ is non-formal.
 \hfill $\Box$
\end{lemma}

Next, we introduce another obstruction to formality, which we call
G-Massey product, since it is a generalization, in spirit, of the
Massey products. 
This product has the advantage of being simpler for computations
than the multiple Massey products.

\begin{definition}\label{def:GMassey}
 Let $X$ be a manifold of any dimension. 
 Let $a, x_1,x_2,x_3 \in H^2(X)$ be degree $2$ cohomology classes
 satisfying that $a\cup x_i=0$, $i=1,2,3$. We define the G-Massey
 product
 $\la a; x_1,x_2,x_3 \ra$ as the subset
  $$ 
  \la a; x_1,x_2,x_3 \ra  = \{ [ \xi_1\wedge\xi_2\wedge\beta_3
  +\xi_2\wedge\xi_3\wedge\beta_1
  +\xi_3\wedge\xi_1\wedge\beta_2] \, |\, a=[\alpha], x_i=[\beta_i],
  \alpha\wedge \beta_i= d\xi_i, i=1,2,3\}
  $$
  $$
  \subset  H^8(X) \, . \qquad\qquad\qquad\qquad\qquad\qquad\qquad\qquad
  \qquad\qquad\qquad
 $$
 We say that the G-Massey product is trivial if $0\in \la a; x_1,x_2,x_3
 \ra$.
\end{definition}

We must notice that
 \begin{equation}\label{eqn:GM2}
 [ \xi_1\wedge\xi_2\wedge\beta_3
  +\xi_2\wedge\xi_3\wedge\beta_1
  +\xi_3\wedge\xi_1\wedge\beta_2]\in H^8(X)
 \end{equation}
is a well defined cohomology class, for forms $\alpha,\beta_i\in
\Omega^2(X)$ and $\xi_i\in \Omega^3(X)$, with $a=[\alpha]$,
$x_i=[\beta_i]$ and $\alpha\wedge \beta_i= d\xi_i$, $i=1,2,3$. In
fact, we have
 $$
  d(\xi_1\wedge\xi_2\wedge\beta_3
  +\xi_2\wedge\xi_3\wedge\beta_1
  +\xi_3\wedge\xi_1\wedge\beta_2)=\alpha\wedge
 \beta_1\wedge \xi_2\wedge \beta_3 - \xi_1 \wedge \alpha\wedge
 \beta_2 \wedge \beta_3 +
 $$
 $$
 +\alpha\wedge \beta_2\wedge
 \xi_3\wedge \beta_1 - \xi_2 \wedge \alpha\wedge \beta_3 \wedge
 \beta_1 +\alpha\wedge \beta_3\wedge \xi_1\wedge \beta_2 - \xi_3
 \wedge \alpha\wedge \beta_1 \wedge \beta_2 =0 \, .
 $$

\begin{lemma}\label{lem:GMassey}
 Let $a, x_1,x_2,x_3 \in H^2(X)$ be degree $2$ cohomology classes
 satisfying that $a\cup x_i=0$, $i=1,2,3$. Then
  $$
 b_1,b_2\in \la a; x_1,x_2,x_3 \ra \Longrightarrow b_1-b_2 \in W ,
 $$
 where $W=\la x_1,a,x_2 \ra \cup H^3(X) +\la x_1,a,x_3 \ra\cup H^3(X) + \la
 x_2,a,x_3 \ra\cup H^3(X) \subset H^8(X)$.
\end{lemma}

\begin{proof}
Choose forms $\alpha,\beta_i\in \Omega^2(X)$ and $\xi_i\in
\Omega^3(X)$, with $a=[\alpha]$, $x_i=[\beta_i]$ and $\alpha\wedge
\beta_i= d\xi_i$, $i=1,2,3$. First of all, note that the
conditions $a\cup x_i=0$, $i=1,2,3$, ensure that the triple Massey
products $\la x_1,a,x_2 \ra$, $\la x_1,a,x_3 \ra$, $\la x_2,a,x_3
\ra$ are well defined.

Now suppose that we write $a=[\alpha+d f]$, $f\in \Omega^1(X)$.
Then $(\alpha+df)\wedge \beta_i= d(\xi_i+ f\wedge \beta_i)$ and
  $$
  (\xi_1+f\wedge \beta_1)\wedge(\xi_2+f\wedge \beta_2)\wedge\beta_3
  +(\xi_2+f\wedge \beta_2)\wedge(\xi_3+f\wedge \beta_3)\wedge\beta_1
  +(\xi_3+f\wedge \beta_3)\wedge(\xi_1+f\wedge \beta_1)\wedge\beta_2=
  $$
  $$
  =\xi_1\wedge\xi_2\wedge\beta_3
  +\xi_2\wedge\xi_3\wedge\beta_1
  +\xi_3\wedge\xi_1\wedge\beta_2,
  $$
so the cohomology class (\ref{eqn:GM2}) does not change by
changing the representative of $a$. If we change the
representatives of $x_i$, say for instance $x_1=[\beta_1+ df]$,
$f\in \Omega^1(X)$, then $\alpha\wedge (\beta_1+df)= d(\xi_1+
\alpha\wedge f)$ and
 $$
  (\xi_1+ \alpha \wedge f)\wedge\xi_2\wedge\beta_3
  +\xi_2\wedge\xi_3\wedge(\beta_1+df)
  +\xi_3\wedge(\xi_1 +\alpha\wedge f)\wedge\beta_2=
  $$
  $$
  =\xi_1\wedge\xi_2\wedge\beta_3
  +\xi_2\wedge\xi_3\wedge\beta_1
  +\xi_3\wedge\xi_1\wedge\beta_2 + d(f\wedge\xi_2\wedge \xi_3),
  $$
  thus the cohomology class (\ref{eqn:GM2}) does not change again.
Finally, if we change the form $\xi_1$ to $\xi_1+ g$, $g\in
\Omega^3(X)$ closed, then
  $$
  (\xi_1+g)\wedge\xi_2\wedge\beta_3
  +\xi_2\wedge\xi_3\wedge\beta_1
  +\xi_3\wedge(\xi_1+g)\wedge\beta_2
  $$
  $$
  =\xi_1\wedge\xi_2\wedge\beta_3
  +\xi_2\wedge\xi_3\wedge\beta_1
  +\xi_3\wedge\xi_1\wedge\beta_2 + g \wedge
  ( \xi_2\wedge\beta_3-\xi_3 \wedge\beta_2),
  $$
and $\xi_2\wedge\beta_3-\xi_3\wedge\beta_2$ is a representative of
$\la x_2,a,x_3 \ra$.
\end{proof}

The indeterminacy of a subset $S$ of a vector space $V$ is the
subspace of $W\subset V$ generated by the differences $s_1-s_2\in
S$ for all $s_1,s_2\in S$. In this situation, $S$ defines an
element in $V/W$. Lemma \ref{lem:GMassey} says that the
indeterminacy of $\la a; x_1,x_2,x_3 \ra$ is contained in $W=\la
x_1,a,x_2 \ra \cup H^3(X) +\la x_1,a,x_3 \ra\cup H^3(X) + \la
x_2,a,x_3 \ra\cup H^3(X)$, hence the G-Massey product $\la a;
x_1,x_2,x_3 \ra$ gives a well-defined element in
  $$
  \frac{H^8(X)}{\la
  x_1,a,x_2 \ra \cup H^3(X) +\la x_1,a,x_3 \ra\cup H^3(X) +
  \la x_2,a,x_3 \ra\cup H^3(X)} \, .
  $$

The relevance of the G-Massey product for formality is given in
the following result.

\begin{proposition}\label{prop:nonformal}
Let $a, x_1,x_2,x_3 \in H^2(X)$ be cohomology classes satisfying
that $a\cup x_i=0$, $i=1,2,3$. Suppose that $\la a; x_1,x_2,x_3
\ra$ is a non-trivial G-Massey product. Then $X$ is not formal.
\end{proposition}

\begin{proof}
Let $\psi:(\bigwedge V,d) \to (\Omega^*(X),d)$ be the minimal
model for $X$. Then there are closed elements $\hat{a}, \hat{x}_i
\in (\bigwedge V)^2$ whose images are $2$--forms $\alpha, \beta_i$
representing $a,x_i$. Since $[\hat{a}\cdot \hat{x}_i]=0$, there
are elements $\hat{\xi}_i\in (\bigwedge V)^3$ such that
$d\hat{\xi}_i= \hat{a}\cdot \hat{x}_i$. Let
$\xi_i=\psi(\hat{\xi}_i)\in \Omega^3(X)$.

If $X$ is formal, then there exists a quasi-isomorphism
$\psi':(\bigwedge V,d) \to (H^*(X),0)$. Note that by adding a
closed element to $\hat{\xi}_i$ we can suppose that
$\psi'(\hat{\xi}_i)=0$. Then
  $$
  [\xi_1\wedge\xi_2\wedge\beta_3
  +\xi_2\wedge\xi_3\wedge\beta_1
  +\xi_3\wedge\xi_1\wedge\beta_1]=\psi'(\hat\xi_1\wedge\hat\xi_2\wedge\hat{x}_3
  +\hat\xi_2\wedge\hat\xi_3\wedge\hat{x}_1
  +\hat\xi_3\wedge\hat\xi_1\wedge\hat{x}_2)=0
  $$
belongs to $\la a; x_1,x_2,x_3 \ra$.
\end{proof}

Actually, the G-Massey product is the first obstruction to
formality that appears as an obstruction to $3$--formality
\cite{FM1} for a simply connected manifold, and which is different
from a Massey product.

The G-Massey product can be related, in some situations, with the
multiple Massey products. However, it cannot be written in terms
of the higher Massey products \cite{K} (or even the matric Massey
products \cite{May,BTe}) because the indeterminacy of Massey
products is usually much bigger. (A similar phenomenon happens to
the product $\la a\ra^k$ discussed in \cite[Section 3]{K}.)

\begin{lemma}\label{lem:quadruple}
Let $X$ be a manifold of any dimension. Let $a, x_1,x_2,x_3 \in
H^2(X)$ be cohomology classes satisfying that $a\cup x_i=0$,
$i=1,2,3$, and $a\cup a=0$. Suppose that $H^5(X)=0$. Then
 \begin{equation}\label{eqn:left}
 \la a; x_1,x_2,x_3\ra \bigcap \Big( x_3  \la x_1,a,a,x_2 \ra +
 x_2 \la x_3,a,a,x_1\ra + x_1 \la x_2,a,a,x_3\ra\Big) \neq \emptyset\,
 .
 \end{equation}
\end{lemma}

\begin{proof}
Write $a=[\alpha]$, $x_i=[\beta_i]$, $\alpha\wedge \beta_i=
d\xi_i$, $i=1,2,3$ and $\alpha\wedge \alpha =d\chi$. The triple
Massey products $\la x_i,a,a\ra$ and $\la a,a,x_i\ra$ are defined
and zero, since $H^5(X)=0$. Then we can write $\xi_i\wedge\alpha -
\beta_i\wedge\chi = d\eta_i$. By definition,
 $$
 \eta_1\wedge \beta_2 -\eta_2\wedge\beta_1 +\xi_1\wedge\xi_2 \in \la x_1,a,a,x_2
 \ra,
 $$
and analogously for the others. Thus the element
 \begin{eqnarray*}
  &&\xi_1\wedge\xi_2\wedge\beta_3
  +\xi_2\wedge\xi_3\wedge\beta_1
  +\xi_3\wedge\xi_1\wedge\beta_2 = \\
 && =\beta_3\wedge(\eta_1\wedge \beta_2 -\eta_2\wedge\beta_1
 +\xi_1\wedge\xi_2)+
 \beta_2\wedge(\eta_3\wedge \beta_1 -\eta_1\wedge\beta_3
 +\xi_3\wedge\xi_1)+ \\
  & & \qquad +\beta_1\wedge(\eta_2\wedge \beta_3 -\eta_3\wedge\beta_2
 +\xi_2\wedge\xi_3)
 \end{eqnarray*}
is in the intersection (\ref{eqn:left}).
\end{proof}

\begin{remark} \label{rem:quadruple}
Note that
 $$
 S= x_3  \la x_1,a,a,x_2 \ra +
 x_2 \la x_3,a,a,x_1\ra + x_1 \la x_2,a,a,x_3\ra \subset
 H^8(X)
 $$
is only defined it $a\cup a=0$. Therefore the G-Massey product can
be understood as a refinement of the subset $S$. Moreover, the
indetermination of $S$ is different (and usually bigger) than that
of the G-Massey product $\la a; x_1,x_2,x_3\ra$.
\end{remark}

The concept of formality is also defined for nilpotent
CW-complexes, and all the discussion above can be extended to them
by using piecewise polynomial differential forms instead of
differential forms. 
Also in this case G-Massey products can be defined.
We shall not need this in full generality, but we shall
use the case when $X$ is an orbifold.

An orbifold is a topological space $X$ with an atlas with charts
modelled on $U/\Pi_p$, where $U$ is an open set of $\RR^n$ and
$\Pi_p$ is a finite group acting linearly on $U$ with only one
fixed point $p\in U$. For an orbifold $X$, we define $\Omega^k(X)$
as the space of orbifold differential forms, i.e., forms such that
in each chart are $\Pi_p$-invariant elements of $\Omega^k(U)$. The
orbifold minimal model of $X$ is defined as the minimal model
$(\bV,d)$ of $(\Omega^k(X),d)$.

\begin{lemma} \label{lem:orb-minimal-model}
Suppose that $X$ is a smooth manifold with minimal model
$(\bV,d)$. Let  $\Pi$ be a finite group acting on $X$ with only
isolated points with non-trivial isotropy, and consider the
orbifold $\widehat X=X/\Pi$.  Let $(\bW,d)$ be the minimal model
of the differential algebra $((\bV)^\Pi, d)$. Then
\begin{itemize}
 \item $(\bW,d)$ is the orbifold minimal model of  $\widehat X$.
 \item Consider $\widehat X$ as a topological space (actually it
 is naturally a CW-complex). If $\widehat X$ is nilpotent, then
 $(\bW,d)$ is its minimal model.
\end{itemize}
\end{lemma}

\begin{proof}
In this situation, $\Omega^k(\widehat X)=\Omega^k(X)^{\Pi}$.
 The action of $\Pi$ on $(\Omega(X),d)$ lifts to an action on the
 minimal model $(\bV,d)$ (see \cite{FelT} for example). As
 $(\bV,d)\to (\Omega(X),d)$ is a quasi-isomorphism,
 $(\bV^\Pi,d)\to (\Omega(X)^\Pi,d)$ is also. Thus $(\bW,d)$ is the
 orbifold minimal model of $X$.

For the second item, triangulate $\widehat{X}$ in such a way that
the orbifolds points are vertices of the triangulation. The
algebra of piecewise polynomial differential forms is
quasi-isomorphic to the algebra $(\Omega^*_{PS}(\widehat X),d)$ of
piecewise smooth differential forms \cite{GM}. Now the natural map
$(\Omega^k(\widehat X),d)\to (\Omega^*_{PS}(\widehat X),d)$ is a
quasi-isomorphism since $H^*(\Omega^k(\widehat X),d) \cong
H^*(X)^{\Pi} \cong H^*(\widehat X)$. So the minimal model of
$\widehat X$ is also $(\bW,d)$.
\end{proof}

\medskip
\section{A nilmanifold of dimension $6$} \label{sec:preliminaries}

Let $G$ be the simply connected nilpotent Lie group of dimension
$6$ defined by the structure equations
 \begin{equation}\label{eqn:a}
 \begin{array}{lll}
  && d\beta_i=0, \qquad i=1,2 \\
  && d\gamma_i=0,\qquad i=1,2 \\
  && d\eta_1=-\beta_1\wedge \gamma_1+ \beta_2\wedge \gamma_1+
     \beta_1\wedge \gamma_2+ 2  \beta_2\wedge \gamma_2, \quad \\
  && d\eta_2=2 \beta_1\wedge \gamma_1+ \beta_2\wedge \gamma_1+
     \beta_1\wedge \gamma_2 - \beta_2\wedge \gamma_2,
 \end{array}
 \end{equation}
where $\{\beta_i,\gamma_i,\eta_i; 1\leq i \leq 2\}$ is a basis of
the left invariant $1$--forms on $G$. Because the structure
constants are rational numbers, Mal'cev theorem \cite{Malc}
implies the existence of a discrete subgroup $\Gamma$ of $G$ such
that the quotient space $N=\Gamma{\backslash} G$ is compact.

Using Nomizu theorem \cite{No} we can compute the real cohomology
of $N$. We get
\begin{eqnarray*}
 H^0(N)& =& \la 1\ra,\\
 H^1(N) &=& \la [\beta_1], [\beta_2],[\gamma_1],[\gamma_2]\ra,\\
 H^2(N) &=& \la [\beta_1 \wedge \beta_2], [\beta_1 \wedge
 \gamma_1],
  [\beta_1 \wedge \gamma_2], [\gamma_1\wedge \gamma_2],
 [\beta_1\wedge \eta_2 - \beta_2\wedge \eta_1],
  [\gamma_1\wedge \eta_2 - \gamma_2\wedge \eta_1],\\
  &&
 [\beta_1\wedge \eta_1 + \beta_1\wedge \eta_2+\beta_2\wedge \eta_2],
  [\gamma_1\wedge \eta_1 + \gamma_1\wedge \eta_2+\gamma_2\wedge \eta_2]\ra,\\
  H^3(N) &=& \la [\beta_1\wedge\beta_2\wedge \eta_1],
 [\beta_1\wedge\beta_2\wedge \eta_2],
 [\gamma_1\wedge\gamma_2\wedge \eta_1],
 [\gamma_1\wedge\gamma_2\wedge \eta_2],
 [\beta_1\wedge\gamma_1\wedge (\eta_1+2\eta_2)],\\
 & &
   [\beta_1 \wedge \gamma_1 \wedge \eta_2-\beta_1 \wedge \gamma_2 \wedge \eta_1],
   [\beta_1 \wedge \gamma_2 \wedge \eta_1-\beta_1 \wedge \gamma_2 \wedge \eta_2],
 [\beta_2 \wedge \gamma_2 \wedge (\eta_2+2\eta_1)],  \\
 & &
  [\beta_2 \wedge \gamma_2 \wedge \eta_1 -\beta_2 \wedge \gamma_1 \wedge \eta_2],
  [\beta_2 \wedge \gamma_1 \wedge \eta_2-\beta_2 \wedge \gamma_1 \wedge   \eta_1]\ra,\\
 H^4(N) &=& \la [\beta_1 \wedge \beta_2 \wedge \gamma_1 \wedge \eta_1],
   [\beta_1 \wedge \beta_2 \wedge \gamma_1 \wedge \eta_2],
  [\beta_1 \wedge \beta_2 \wedge \eta_1 \wedge \eta_2],
 [\beta_1 \wedge \gamma_1 \wedge \gamma_2 \wedge \eta_2],\\
  & &   [\beta_2 \wedge \gamma_1 \wedge \gamma_2 \wedge \eta_2],
   [\gamma_1 \wedge \gamma_2 \wedge \eta_1 \wedge \eta_2],
   [\beta_1 \wedge \gamma_2 \wedge \eta_1 \wedge \eta_2-
   \beta_2 \wedge \gamma_1 \wedge \eta_1 \wedge \eta_2],\\
  & & [\beta_1 \wedge \gamma_2 \wedge \eta_1 \wedge \eta_2+
 \beta_1 \wedge \gamma_1 \wedge \eta_1 \wedge \eta_2+
 \beta_2 \wedge \gamma_2 \wedge \eta_1 \wedge \eta_2]\ra,\\
 H^5(N) &=& \la [\beta_1 \wedge \beta_2 \wedge \gamma_1 \wedge
  \eta_1 \wedge \eta_2],
 [\beta_1 \wedge \beta_2 \wedge \gamma_2 \wedge
  \eta_1 \wedge \eta_2],
  [\beta_1 \wedge \gamma_1 \wedge \gamma_2 \wedge
  \eta_1 \wedge \eta_2],\\
 & &
  [\beta_2 \wedge \gamma_1 \wedge \gamma_2 \wedge
  \eta_1 \wedge \eta_2] \ra, \\
  H^6(N) &=& \la [\beta_1 \wedge \beta_2 \wedge \gamma_1 \wedge \gamma_2
  \wedge \eta_1 \wedge \eta_2]\ra.
 \end{eqnarray*}

We can give a more explicit description of the group $G$. As a
differentiable manifold $G=\RR^6$. The nilpotent Lie group
structure of $G$ is given by the multiplication law
 \begin{equation}\label{eqn:m}
 \begin{array}{ccl}
  m:\qquad G \times G&\too& \qquad G \\
  \left((y_1',y_2',z_1',z_2',v_1',v_2'),(y_1,y_2,z_1,z_2,v_1,v_2)\right)
  &\mapsto & \Big(y_1+y_1',y_2+y_2',z_1+z_1',z_2+z_2', \\ &&
  v_1 +v_1'+(y_1' -y_2') z_1 - (y_1'+2y_2') z_2,\\ &&
  v_2 +v_2'- (2y_1' +y_2') z_1 +(y_2'-y_1') z_2 \Big).
 \end{array}
 \end{equation}

We also need a discrete subgroup, which it could be taken to be
$\ZZ^6 \subset G$. However, for later convenience, we shall take
the subgroup
 $$
 \Gamma = \{(y_1,y_2,z_1,z_2,v_1,v_2)\in \ZZ^6 \, |\, v_1 \equiv v_2
 \pmod 3 \} \subset G,
 $$
and define the nilmanifold
 $$
 N=\Gamma \backslash G\ .
 $$
In terms of a (global) system of coordinates
$(y_1,y_2,z_1,z_2,v_1,v_2)$ for $G$, the $1$--forms $\beta_i$,
$\gamma_i$ and $\eta_i$, $1\leq i \leq 2$, are  given by
 \begin{eqnarray*}
  \beta_i&=&dy_i, \quad  1\leq i \leq 2,  \quad \\
  \gamma_i &= &dz_i,  \quad  1\leq i \leq 2,  \quad \\
  \eta_1&=&dv_1-y_1 dz_1 +y_2dz_1+y_1 dz_2 +2y_2dz_2, \\
  \eta_2&=&dv_2+2y_1 dz_1 +y_2dz_1+y_1 dz_2 -y_2dz_2.
 \end{eqnarray*}
Note that $N$ is a principal torus bundle
 $$
 T^2=\ZZ\la (1,1),(3,0)\ra
 \backslash \RR^2 \inc N \longrightarrow
 T^4= \ZZ^4\backslash \RR^4,
 $$
with the projection $(y_1,y_2,z_1,z_2,v_1,v_2) \mapsto
(y_1,y_2,z_1,z_2)$.

\medskip

The Lie group $G$ can be also described as follows. Consider the
basis $\{\mu_i,\nu_i,\theta_i; 1\leq i \leq 2\}$ of the left
invariant $1$--forms on $G$ given by
 $$
 \begin{array}{ll}
 \mu_1 = \beta_1+\displaystyle\frac{1+\sqrt{3}}{2} \beta_2,  \quad &
 \mu_2 = \beta_1+ \displaystyle\frac{1-\sqrt{3}}{2} \beta_2,  \quad \\
 \nu_1 = \gamma_1+\displaystyle\frac{1+\sqrt{3}}{2} \gamma_2,  \quad &
 \nu_2 = \gamma_1+ \displaystyle\frac{1-\sqrt{3}}{2} \gamma_2,  \quad \\
 \theta_1 = \displaystyle\frac{2}{\sqrt{3}} \eta_1+
 \displaystyle\frac{1}{\sqrt{3}} \eta_2, \quad  & \theta_2=\eta_2.
 \end{array}
 $$
Hence, the structure equations can be rewritten as
 \begin{equation} \label{eqn:struc}
 \begin{array}{l}
  d\mu_i=0,  \quad  1\leq i \leq 2,  \\
  d\nu_i=0,  \quad  1\leq i \leq 2,  \\
  d\theta_1=\mu_1\wedge \nu_1- \mu_2\wedge \nu_2,  \\
  d\theta_2=\mu_1\wedge \nu_2+ \mu_2\wedge \nu_1.
 \end{array}
 \end{equation}
This means that $G$ is the complex Heisenberg group $H_{\CC}$,
that is, the complex nilpotent Lie group of complex matrices of
the form
 $$
 \pmatrix{1&u_2&u_3\cr 0&1&u_1\cr 0&0&1\cr}.
 $$
In fact, in terms of the natural (complex) coordinate functions
$(u_1,u_2,u_3)$ on $H_{\CC}$, we have that the complex $1$--forms
 $$
 \mu=du_1, \ \nu=du_2, \ \theta=du_3-u_2 du_1
 $$
are left invariant and $d\mu=d\nu=0$, $d\theta=\mu\wedge\nu$. Now,
it is enough to take $\mu_1=\Re(\mu)$, $\mu_2=\Im(\mu)$,
$\nu_1=\Re(\nu)$, $\nu_2=\Im(\nu)$, $\theta_1=\Re(\theta)$,
$\theta_2=\Im(\theta)$ to recover equations (\ref{eqn:struc}).

\begin{lemma} \label{lem:N}
  Let $\Lambda \subset \CC$ be the lattice generated by $1$ and
  $\zeta=e^{2\pi i/3}$, and consider the discrete subgroup
  $\Gamma_H \subset H_{\CC}$ formed
  by the matrices in which $u_1,u_2,u_3 \in \Lambda$. Then
  there is a natural identification of
  $N=\Gamma{\backslash}G$ with the quotient $\Gamma_H\backslash
  H_{\CC}$.
\end{lemma}

\begin{proof}
  We have constructed above an isomorphism of Lie groups $G\to H_{\CC}$, whose
  explicit equations are
   $$
   (y_1,y_2,z_1,z_2,v_1,v_2) \mapsto (u_1,u_2,u_3),
   $$
where
 \begin{eqnarray*}
 u_1 &=& \left( y_1+\displaystyle\frac{1+\sqrt{3}}{2} y_2  \right) +
 i \left( y_1+ \displaystyle\frac{1-\sqrt{3}}{2} y_2 \right), \\
 u_2 &=& \left( z_1+\displaystyle\frac{1+\sqrt{3}}{2} z_2  \right) +
 i \left( z_1+ \displaystyle\frac{1-\sqrt{3}}{2} z_2 \right), \\
 u_3 &=& \frac{1}{\sqrt{3}}\left( 2v_1+v_2 +3z_1y_2 +3z_2y_1+3z_2y_2
 \right) + i \left(v_2 + 2z_1y_1+z_2y_1 + z_1y_2 -z_2y_2\right).
 \end{eqnarray*}
Note that the formula for $u_3$ can be deduced from
 $$
 du_3-u_2du_1 =\theta =\left( \displaystyle\frac{2}{\sqrt{3}}\eta_1 +
 \displaystyle\frac{1}{\sqrt{3}}\eta_2 \right) + i\eta_2 \ .
 $$
Now the group $\Gamma\subset G$ corresponds under this
isomorphism to
 $$
 \left\{(u_1,u_2,u_3) | u_1, u_2 \in \ZZ\left\la 1+i,
  \frac{1+\sqrt{3}}{2}+ \frac{1-\sqrt{3}}{2} i \right\ra
  , u_3 \in \ZZ\left\la
  2\sqrt{3}, \sqrt{3} + i \right\ra \right\}.
  $$
Using the isomorphism of Lie groups  $H_\CC \to H_\CC$ given by
 $$
 (u_1,u_2,u_3) \mapsto (u_1',u_2',u_3')=\left( \frac{u_1}{1+i},\frac{u_2}{1+i}
  ,\frac{u_3}{(1+i)^2}\right),
  $$
we get that $u_1',u_2', u_3'\in \Lambda=\ZZ\la 1,\zeta\ra$, which
completes the proof.
\end{proof}

\begin{remark}
 If we had considered the discrete subgroup $\ZZ^6\subset G$ instead of
 $\Gamma\subset G$, then we would not have obtained the fact $u_3'\in
 \Lambda$ in the proof of Lemma \ref{lem:N}. Actually the
 manifold $\ZZ^6\backslash G$ is not diffeomorphic to
 $N=\Gamma\backslash G$, as can be proved as in Proposition
 \ref{prop:iwasawa} below. Note that $N=\Gamma\backslash G \surj
 \ZZ^6\backslash G$ is a $3:1$ covering.
 (However the nilmanifold $\ZZ^6\backslash G$ could
 also have been used as a starting point to construct a
 simply connected compact symplectic non-formal $8$--manifold with
 the arguments of \cite{FM4}.)
\end{remark}

Under the identification $N=\Gamma{\backslash}G \cong
\Gamma_H\backslash H_{\CC}$, $N$ is a principal torus bundle
 $$
 T^2={\Lambda \backslash \CC} \inc N \longrightarrow
 T^4= {\Lambda^2 \backslash \CC^2},
 $$
via the projection $(u_1,u_2,u_3)\mapsto (u_1,u_2)$.

The manifold $(\Gamma_H\backslash H_{\CC}) \times (\Lambda
\backslash\CC)$ is the $8$--dimensional compact nilmanifold $M$
defined in \cite[Section 2]{FM4}. It is interesting here to
compare $N$ with the Iwasawa manifold. Let us recall its
definition. Let $\Lambda'\subset \CC$ be the Gaussian integers,
i.e., the lattice generated by $1$ and $i$, and consider the
discrete subgroup $\Gamma_0 \subset H_{\CC}$ formed by the
matrices in which $u_1,u_2,u_3 \in \Lambda'$. Then the Iwasawa
manifold is defined as the quotient \cite{Chern, GHarris, MK}
 $$
 N'=\Gamma_0 \backslash H_{\CC}.
 $$

Note that $N{}'$ is also a principal torus bundle
 $$
 T^2={\Lambda{}' \backslash \CC} \inc N' \longrightarrow
 T^4=({\Lambda{}')^2 \backslash \CC^2},
 $$

\begin{proposition} \label{prop:iwasawa}
The fundamental groups $\pi_1(N)$ and $\pi_1(N{}')$
are not isomorphic. In particular, $N$ and $N{}'$ are not
diffeomorphic.
\end{proposition}

\begin{proof}
Let us first consider the manifold $N$. It is a principal torus
bundle over $T^4=\Lambda^2 \backslash \CC^2$, the action of
$T^2=\Lambda\backslash \CC$ being by translations in the $u_3$
coordinate. The $1$--form $\theta =du_3 -u_2 du_1 \in
\Omega^1(N,\CC)$ is a connection $1$--form with values in $\CC$,
the Lie algebra of $T^2$. The curvature form $F=d\theta$ is the
lift of the $2$--form $\mu \wedge \nu =d u_1 \wedge du_2 \in
\Omega^2(T^4,\CC^2)$. The cohomology class defined by the
curvature is
 $$
 [F]\in H^2(T^4,\CC)=\Hom (H_2(T^4,\ZZ), \CC).
 $$
The image of this map lies in $\Lambda\subset \CC$. Actually, the
$T^2=\Lambda \backslash \CC$-principal bundles over a space $X$ are
classified by
 $$
 [X, B( \Lambda\backslash \CC)]= [X, K(\Lambda,2)]=H^2(X,\Lambda),
 $$
and $[F]$ gives the required element classifying $N$. Since
$H_k(T^4,\ZZ) = \bigwedge^k (\Lambda^2)$, we can view
intrinsically
 $$
 [F] \in \Hom (\bigwedge\nolimits^2 (\Lambda^2), \Lambda).
 $$

To compute $[F]$, consider the basis for $H_1(T^4,\ZZ)$ given as
$\{e_1=(1,0),e_2=(\zeta,0),e_3=(0,1),e_4=(0,\zeta)\}$. This basis
gives us an isomorphism $H_1(T^4,\ZZ)=\Lambda^2 \iso \ZZ^4$. Then
a basis for the $2$--homology $H_2(T^4,\ZZ)\iso \bigwedge^2 \ZZ^4$
is $\{e_1\wedge e_2,e_1\wedge e_3,e_1\wedge e_4,e_2\wedge
e_3,e_2\wedge e_4, e_3\wedge e_4\}$. Also use the basis
$\{1,\zeta\}$ for $\Lambda$. We compute $[F] \in H^2(T^4,
\Lambda)=\Hom (\bigwedge^2 (\Lambda^2), \Lambda)\cong \Hom
(\bigwedge^2 (\ZZ^4), \ZZ)$ in terms of these bases:
\begin{eqnarray*}
 \ [F]( e_1\wedge e_2) &=& \int_{e_1\wedge e_2} F =
  \int_{u_1\in ({\Lambda \backslash \CC}),\ u_2=0} d u_1\wedge du_2=0, \\
 \ [F]( e_1\wedge e_3) &=& \int_{e_1\wedge e_3} F =
  \int_{u_1=t_1, u_2=t_2, \ 0\leq t_1,t_2\leq 1} du_1\wedge du_2
  = \int_{0}^{1}\int_{0}^{1} d t_1 dt_2=1, \\
 \ [F]( e_1\wedge e_4) &=& \int_{e_1\wedge e_4} F =
  \int_{u_1=t_1, u_2=t_2\zeta,\ 0\leq t_1,t_2\leq 1}  d u_1\wedge du_2
  = \int_{0}^{1}\int_{0}^{1} \zeta d t_1 dt_2=\zeta, \\
 \ [F]( e_2\wedge e_3) &=& \int_{e_2\wedge e_3} F =
  \int_{u_1=t_1\zeta, u_2=t_2,\ 0\leq t_1,t_2\leq 1}  d u_1\wedge du_2
  = \int_{0}^{1}\int_{0}^{1} \zeta d t_1 dt_2=\zeta, \\
 \ [F]( e_2\wedge e_4) &=& \int_{e_2\wedge e_4} F =
  \int_{u_1=t_1\zeta, u_2=t_2\zeta,\ 0\leq t_1,t_2\leq 1}  d u_1\wedge du_2
  = \int_{0}^{1}\int_{0}^{1} \zeta^2 d t_1 dt_2=\zeta^2=-1-\zeta, \\
 \ [F]( e_3\wedge e_4) &=& \int_{e_3\wedge e_4} F =
 \int_{u_1=0,\ u_2\in ({\Lambda \backslash \CC})} d u_1\wedge
 du_2=0.
\end{eqnarray*}
In terms of the given bases, $[F]$ is the matrix
 \begin{equation}\label{eqn:F}
 [F]= \pmatrix{0&1&0&0&-1&0\cr 0&0&1&1&-1&0\cr}.
 \end{equation}

We can similarly work out the case of the Iwasawa manifold $N'$.
Again it is a principal $T^2$-torus bundle over $T^4$, where
$T^2=\Lambda'\backslash \CC$ and $T^4=(\Lambda')^2\backslash
\CC^2$. Working analogously as before, the curvature $F'$ of this
principal bundle is $F'=du_1\wedge du_2$ and the cohomology class
$[F']\in \mathrm{Hom}\,(\bigwedge^2\left((\Lambda')^2\right),
\Lambda')$ is computed as follows: consider the basis
$\{e_1=(1,0),e_2=(i,0),e_3=(0,1),e_4=(0,i)\}$ for $H_1(T^4,\ZZ)$
and the basis $\{1,i\}$ for $\Lambda'$. Then
\begin{eqnarray*}
 \ [F']( e_1\wedge e_2) &=& \int_{e_1\wedge e_2} F =
  \int_{u_1\in ({\Lambda' \backslash \CC}),\ u_2=0} d u_1\wedge du_2=0, \\
 \ [F']( e_1\wedge e_3) &=& \int_{e_1\wedge e_3} F =
  \int_{u_1=t_1, u_2=t_2, \ 0\leq t_1,t_2\leq 1} du_1\wedge du_2
  = \int_{0}^{1}\int_{0}^{1} d t_1 dt_2=1, \\
 \ [F']( e_1\wedge e_4) &=& \int_{e_1\wedge e_4} F =
  \int_{u_1=t_1, u_2=t_2 i,\ 0\leq t_1,t_2\leq 1}  d u_1\wedge du_2
  = \int_{0}^{1}\int_{0}^{1} i d t_1 dt_2=i, \\
 \ [F']( e_2\wedge e_3) &=& \int_{e_2\wedge e_3} F =
  \int_{u_1=t_1i, u_2=t_2,\ 0\leq t_1,t_2\leq 1}  d u_1\wedge du_2
  = \int_{0}^{1}\int_{0}^{1} i d t_1 dt_2=i, \\
 \ [F']( e_2\wedge e_4) &=& \int_{e_2\wedge e_4} F =
  \int_{u_1=t_1i, u_2=t_2i,\ 0\leq t_1,t_2\leq 1}  d u_1\wedge du_2
  = \int_{0}^{1}\int_{0}^{1} i^2 d t_1 dt_2=-1, \\
 \ [F']( e_3\wedge e_4) &=& \int_{e_3\wedge e_4} F =
 \int_{u_1=0,\ u_2\in ({\Lambda' \backslash \CC})} d u_1\wedge
 du_2=0.
\end{eqnarray*}
So the corresponding matrix is
 \begin{equation}\label{eqn:F'}
 [F']= \pmatrix{0&1&0&0&-1&0\cr 0&0&1&1&0&0\cr}.
 \end{equation}

Since (\ref{eqn:F}) and (\ref{eqn:F'}) are different, the two
torus bundles are not isomorphic (as principal bundles over
$T^4$). Moreover they are not even isomorphic even after an
automorphism of the basis, since $[F]$ and $[F{}']$ are
inequivalent under the natural action of $\GL(4,\ZZ) \times
\GL(2,\ZZ)$ (by changes of basis in the lattices corresponding to
base and fiber, respectively). This can be seen by considering the
intersection pairing $Q$ on $\bigwedge^2 \ZZ^4 \times \bigwedge^2
\ZZ^4 \to \bigwedge^4 \ZZ^4 \cong \ZZ$ (well-defined modulo sign)
and looking at the determinant of the image lattice of $[F]$ and
$[F']$, respectively. As $\Im
[F]=\{(0,1,0,0,-1,0),(0,0,1,1,-1,0)\}$ and $\Im
[F{}']=\{(0,1,0,0,-1,0),(0,0,1,1,0,0)\}$, then
 $$
 \det(Q\vert_{\Im [F]})=\det\pmatrix{2&1\cr 1&2\cr}=3
 $$
and
 $$
 \det(Q\vert_{\Im [F{}']})=\det\pmatrix{2&0\cr 0&2\cr}=4.
 $$

A relevant point here is that the fundamental group can be read
off from the classifying cohomology class. For instance, the
fundamental group of $N$ is an extension
 $$
 \Lambda =\pi_1(T^2) \to \Gamma_H=\pi_1(N) \to \Lambda^2 =\pi_1(T^4)
 $$
and this is determined by the commutator bracket
 $$
 [\ ,\ ]: \Lambda^2 \times \Lambda^2  \too \Lambda,
 $$
which in turn coincides with the linear map $[F]$. Note that for
$N$, $\pi_1(T^2)$ is exactly the center of $\pi_1(N)$. This holds
since $[F]=[\ ,\ ]$ is non-degenerate. This means that
$\pi_1(T^2)$ and $\pi_1(T^4)$ are univocally determined by
$\pi_1(N)$. Hence the orbit of $[F]$ under $\GL(4,\ZZ) \times
\GL(2,\ZZ)$ determines the fundamental group $\pi_1(N)$. A similar
fact happens for $\pi_1(N')=\Gamma_0$. Therefore, $\pi_1(N)
\not\cong \pi_1(N{}')$, and $N$ and $N{}'$ are not diffeomorphic.
\end{proof}

\medskip
\section{Quotient of a nilmanifold by a $\ZZ_3$-action} \label{sec:hatE}

We define the $8$--dimensional compact nilmanifold $M$ as  the
product
 $$
  M=T^2 \times N.
 $$
By Lemma \ref{lem:N} there is an isomorphism between $M$ and the
manifold $(\Gamma_H\backslash H_{\CC}) \times (\Lambda \backslash
\CC)$ studied in \cite[Section 2]{FM4} (we have to send the factor
$T^2$ of $M$ to the factor $\Lambda \backslash \CC$). Clearly, $M$
is a principal torus bundle
 $$
 T^2 \inc M \stackrel{\pi}{\too} T^6.
 $$
Let $(x_1,x_2)$ be the Lie algebra coordinates for $T^2$, so that
$(x_1,x_2,y_1,y_2,z_1,z_2,v_1,v_2)$ are coordinates for the Lie
algebra $\RR^2\times G$ of $M$. Then
$\pi(x_1,x_2,y_1,y_2,z_1,z_2,v_1,v_2)=(x_1,x_2,y_1,y_2,z_1,z_2)$.
A basis for the left invariant (closed) $1$--forms on $T^2$ is
given as $\{\alpha_1,\alpha_2\}$, where $\alpha_1=dx_1$ and
$\alpha_2=dx_2$. Then $\{\alpha_i,\beta_i,\gamma_i,\eta_i; 1\leq i
\leq 2\}$ constitutes a (global) basis for the left invariant
$1$--forms on $M$. Note that $\{\alpha_i,\beta_i,\gamma_i; 1\leq i
\leq 2\}$ is a basis for the left invariant closed $1$--forms on
the base $T^6$. (We use the same notation for the differential
forms on $T^6$ and their pullbacks to $M$.) Using the computation
of the cohomology of $N$, we get that the Betti numbers of $M$
are: $b_0(M)=b_8(M)=1$, $b_1(M)=b_7(M)=6$, $b_2(M)=b_6(M)=17$,
$b_3(M)=b_5(M)=30$, $b_4(M)=36$. In particular, $\chi(M)=0$, as
for any nilmanifold.

\medskip

Let us now write the minimal model of the nilmanifold $M$.
Nomizu's theorem \cite{No} gives that the minimal model of $M$ is
the differential graded commutative algebra
 $$
  (\bigwedge W,d)=(\bigwedge(a_1,a_2, b_1,b_2 ,c_1,c_2,e_1, e_2 ),d),
 $$
whose generators $a_i$, $b_i$, $c_i$ and $e_i$, $1\leq i \leq 2$, have
degree $1$, the differential $d$ is given by
 \begin{eqnarray*}
  && da_i=db_i=dc_i=0,  \quad  1\leq i \leq 2,\\
  && de_1=-b_1 \cdot c_1+ b_2 \cdot c_1+ b_1 \cdot c_2 + 2 b_2 \cdot c_2, \\
  && de_2=2b_1 \cdot c_1+ b_2 \cdot c_1+ b_1 \cdot c_2 - b_2 \cdot c_2,
 \end{eqnarray*}
and the morphism $\phi \colon  (\bigwedge(a_i,b_i,c_i,e_i),d)\to
(\Omega(M),d)$, inducing an isomorphism on cohomology, is defined
by $\phi(a_i)=\alpha_i$, $\phi(b_i)=\beta_i$,
$\phi(c_i)=\gamma_i$, $\phi(e_i)=\eta_i$, for $1\leq i\leq 2$,
where $(\Omega(M),d)$ denotes the de Rham complex of differential
forms on $M$.

\medskip

Consider the action of the finite group $\ZZ_{3}$ on ${\RR^2}$ given by
 $$
  \rho(x_1,x_2)=(-x_1-x_2,x_1),
 $$
for $(x_1,x_2)\in {\RR^2}$, $\rho$ being the generator of $\ZZ_3$.
Clearly $\rho(\ZZ^{2})=\ZZ^{2}$, and so $\rho$ defines an action
of $\ZZ_{3}$ on the $2$-torus $T^2=\ZZ^2\backslash \RR^2$ with $3$
fixed points: $(0,0)$, $(\frac13,\frac13)$ and
$(\frac23,\frac23)$. The quotient space $T^2/\ZZ_{3}$ is the
orbifold $2$--sphere $S^2$ with $3$ points of multiplicity $3$.
Let $x_1$, $x_2$ denote the natural coordinates functions on
${\RR^2}$. Then the $1$--forms $dx_1$, $dx_2$ satisfy
$\rho^*(dx_1)=-dx_1-dx_2$ and $\rho^*(dx_2)=dx_1$, hence
$\rho^*(-dx_1-dx_2)=dx_2$. Thus, we can take the $1$-forms
$\alpha_1$ and $\alpha_2$ on $T^2$ such that
 \begin{equation} \label{eqn:v1}
   \rho^*(\alpha_1)=-\alpha_1-\alpha_2, \quad
   \rho^*(\alpha_2)=\alpha_1.
 \end{equation}

We denote by $A$ the $2$-dimensional representation of $\ZZ_3$ given by
 \begin{equation} \label{eqn:matrix}
  \begin{array}{ccl}
   \ZZ_3 & \too & \GL(2,\RR) \\
   \rho & \mapsto & \left(\begin{array}{cc} -1 &-1 \\ 1&0 \end{array}\right)
  \end{array}
  \end{equation}
Then the cohomology group $H^1(T^2)\iso A$, as
$\ZZ_3$-representations.

It is easy to see the following isomorphisms of representations \cite{FultonHarris}:
 \begin{equation} \label{eqn:iso}
 A\wedge A \iso {\RR},  \qquad A \otimes A \iso  {\RR} \oplus {\RR} \oplus A,
 \end{equation}
where $\RR$ denotes the trivial $1$--dimensional representation.

Define the following action of $\ZZ_3$ on $M$, given, at the level
of Lie groups, by $\rho \colon {\RR^2} \times
{\RR^6}\longrightarrow {\RR^2} \times {\RR^6}$,
 $$
 \rho(x_1,x_2,y_1,y_2,z_1,z_2,v_1,v_2)
 =(-x_1-x_2,x_1,-y_1-y_2,y_1,-z_1-z_2,z_1,-v_1-v_2,v_1).
 $$
Note that $m(\rho(p'),\rho(p)) = \rho (m (p',p))$, for all
$p,p'\in G$, where $m$ is the multiplication map (\ref{eqn:m}) for
$G$. Also $\Gamma\subset G$ is stable by $\rho$ since
 $$
 v_1\equiv v_2 \pmod 3 \Longrightarrow -v_1-v_2\equiv v_1 \pmod 3.
 $$
Therefore there is a induced map $\rho \colon M \to M$, and this
covers the action $\rho: T^6 \to T^6$ on the $6$--torus
$T^{6}=T^{2}\times T^{2}\times T^{2}$ (defined as the action
$\rho$ on each of the three factors simultaneously). The action of
$\rho$ on the fiber $T^{2}= \ZZ\la (1,1),(3,0)\ra$ has also $3$
fixed points: $(0,0)$, $(1,0)$ and $(2,0)$. Hence there are
$3^4=81$ fixed points on $M$.

\begin{remark}
Under the isomorphism $M \cong (\Gamma_H\backslash H_{\CC}) \times
(\Lambda \backslash \CC)$,
 we have that the action of $\rho$ becomes $\rho(u_1,u_2,u_3)=(\bar \zeta
 u_1,\bar \zeta u_2,\zeta u_3)$, where
 $\zeta=e^{2\pi i/3}$. Composing the isomorphism of Lemma
 \ref{lem:N} with the conjugation $(u_1,u_2,u_3)\mapsto
 (v_1,v_2,v_3)=(\bar{u}_1,\bar{u}_2,
 \bar{u}_3)$ (which is an isomorphism of Lie groups $H_\CC \to H_\CC$ leaving $\Gamma_H$
 invariant),
 we have that the action of $\rho$ becomes
 $\rho(v_1,v_2,v_3)=(\zeta v_1, \zeta v_2,\zeta^2 v_3)$. This is
 the action used in \cite{FM4}.
\end{remark}

We take the basis $\{\alpha_i,\beta_i,\gamma_i,\eta_i; 1\leq i
\leq 2\}$ of the $1$--forms on $M$ considered above. The
$1$--forms $dy_i$, $dz_i$, $dv_i$, $1\leq i\leq 2$, on $G$ satisfy
the following conditions similar to (\ref{eqn:v1}):
$\rho^*(dy_1)=-dy_1-dy_2$, $\rho^*(dy_2)=dy_1$,
$\rho^*(dz_1)=-dz_1-dz_2$, $\rho^*(dz_2)=dz_1$,
$\rho^*(dv_1)=-dv_1-dv_2$, $\rho^*(dv_2)=dv_1$. So
 \begin{equation} \label{eqn:v2}
 \begin{array}{ll}
  \rho^*(\alpha_1)=-\alpha_1-\alpha_2, \quad   & \rho^*(\alpha_2)=\alpha_1,\\
  \rho^*(\beta_1)=-\beta_1-\beta_2, \quad  & \rho^*(\beta_2)=\beta_1,\\
  \rho^*(\gamma_1)=-\gamma_1-\gamma_2,  \quad  &\rho^*(\gamma_2)=\gamma_1,\\
  \rho^*(\eta_1)=-\eta_1-\eta_2,  \quad &\rho^*(\eta_2)=\eta_1.
 \end{array}
 \end{equation}

\begin{remark} \label{rem:3}
If we define the $1$--forms $\alpha_3=-\alpha_1-\alpha_2$,
$\beta_3=-\beta_1-\beta_2$, $\gamma_3=-\gamma_1-\gamma_2$ and
$\eta_3=-\eta_1-\eta_2$, then we have
$\rho^*(\alpha_1)=\alpha_3,\, \rho^*(\alpha_2)=\alpha_1,\,
\rho^*(\alpha_3)=\alpha_2$, and analogously for the others.
\end{remark}

Note that there is also a $\ZZ_3$-action on the minimal model
$(\bigwedge W,d)$ of $M$ defined analogously to (\ref{eqn:v2}). As
$\ZZ_3$-representations, we have an isomorphism $W \iso A^4$. This
gives, using (\ref{eqn:iso}), the following decomposition of the
minimal model as $\ZZ_3$-representation:
  \begin{equation}\label{eqn:model}
  \left\{ \begin{array}{lcl}
   (\bigwedge W)^1 &\iso& A^4,\\
   (\bigwedge W)^2 &\iso& \RR^{16} \oplus A^{6},\\
   (\bigwedge W)^3 &\iso& \RR^{8} \oplus A^{24},\\
   (\bigwedge W)^4 &\iso& \RR^{36} \oplus A^{17},\\
   (\bigwedge W)^5 &\iso& \RR^{8} \oplus A^{24},\\
   (\bigwedge W)^6 &\iso& \RR^{16} \oplus A^{6},\\
   (\bigwedge W)^7 &\iso& A^{4},\\
   (\bigwedge W)^8 &\iso& {\RR}.
  \end{array} \right.
  \end{equation}

\medskip

Define the quotient space
 $$
 \widehat{M}=M/\ZZ_{3},
 $$
and denote by $\varphi:M\to \widehat{M}$ the projection. It is an
orbifold, but we can compute the rational homotopy type of the
underlying topological manifold using Lemma
\ref{lem:orb-minimal-model}. A model for $\widehat{M}$ is given by
the $\ZZ_3$-invariant part $((\bigwedge W)^{\ZZ_3},d)$ of the
minimal model of $M$. This corresponds to the $\RR$-factors of
(\ref{eqn:model}). Since $(\bigwedge W)^1 =W \iso A^4$, the
invariant part $W^{\ZZ_3}$ is zero. This means that the first
stage of the minimal model of $\widehat{M}$ is zero and hence
$b_1(\widehat{M})=0$. This was the starting point that led us to
consider the equations (\ref{eqn:a}) to define $M$.

One can compute explicitly the differential $d:((\bigwedge
W)^i)^{\ZZ_3} \to ((\bigwedge W)^{i+1})^{\ZZ_3}$ to get the
cohomology of $\widehat{M}$. For instance,
 $$
  \begin{array}{lcl}
  H^1(\widehat{M}) &=& 0, \\
  H^2(\widehat{M}) &=& \la [\alpha_1 \wedge \alpha_2],
  [\alpha_1 \wedge \beta_2-\alpha_2 \wedge \beta_1],
    [\alpha_1 \wedge \beta_1+\alpha_1 \wedge \beta_2+\alpha_2 \wedge \beta_2],
    \\
   &&[\alpha_1 \wedge \gamma_2-\alpha_2 \wedge \gamma_1],
    [\alpha_1 \wedge \gamma_1+\alpha_1 \wedge \gamma_2+\alpha_2 \wedge \gamma_2],
   [\beta_1 \wedge \beta_2],
     [\beta_1 \wedge \gamma_2-\beta_2 \wedge \gamma_1],
  \\
 && [\beta_1 \wedge\gamma_1+\beta_1 \wedge\gamma_2+\beta_2 \wedge\gamma_2],
 [\beta_1\wedge \eta_2 - \beta_2\wedge \eta_1],
 [\beta_1\wedge \eta_1 + \beta_1\wedge \eta_2+\beta_2\wedge \eta_2],
 \\
 &&  [\gamma_1\wedge \gamma_2],[\gamma_1\wedge \eta_2 - \gamma_2\wedge \eta_1],
  [\gamma_1\wedge \eta_1 + \gamma_1\wedge \eta_2+\gamma_2\wedge \eta_2]\ra,\\
  H^3(\widehat{M}) &=& 0.
  \end{array}
  $$

\begin{remark}
The Euler characteristic of $\widehat{M}$ can be computed via the
formula for finite group action quotients: let $\Pi$ be the cyclic
group of order $n$, acting on a space $X$ almost freely. Then
 $$
 \chi(X/\Pi)= \frac1n \chi(X) + \sum_{p} \left(1- \frac{1}{\#
 \Pi_p}\right),
 $$
where $\Pi_p\subset \Pi$ is the isotropy group of $p\in X$. In our
case $\chi(\widehat{M})= \frac13 \chi (M) + 81 (1-\frac13) =54$.
\end{remark}

Using this remark and the previous calculation, we get that
$b_1(\widehat{M})=b_7(\widehat{M})=0$, $b_2(\widehat{M}) =
b_6(\widehat{M}) ={13}$, $b_3(\widehat{M}) =b_5(\widehat{M}) =0$
and $b_4(\widehat{M}) =26$. Note that $\widehat{M}$ satisfies
Poincar{\'e} duality since
 $$
 H^*(\widehat{M})= H^*(M)^{\ZZ_3}
 $$
and $H^*(M)$ satisfies Poincar{\'e} duality.

\begin{proposition}\label{prop:thm_simply connected0}
  $\widehat{M}$ is simply connected.
\end{proposition}

\begin{proof}
Let $p_0\in M$ be a fixed point of the $\ZZ_3$-action and let
$\hat{p}_0=\varphi(p_0)$. There is (see \cite{Bredon}) an
epimorphism of fundamental groups
 $$
 \Gamma=\pi_1(M,p_0) \surj \pi_1(\widehat{M},\hat{p}_0).
 $$
This holds since every path in $\widehat{M}$ can be lifted to $M$,
in an unique way as long as it does not touch a singular point, an
in three different ways when it does.

Since the nilmanifold $M$ is a principal torus bundle over the
$6$--torus $T^6$, we have
 $$
 \ZZ^2 \inc \Gamma \to \ZZ^6.
 $$
Consider $p_0 \in M$ a fixed point of the $\ZZ_3$-action and
$\bar{p}_0=\pi(p_0)$, where $\pi \colon M \to T^6$ is the
projection of the torus bundle. Then $\ZZ_3$ acts on
$\pi^{-1}(\bar{p}_0)\iso T^2$, and the restriction to $\ZZ^2$ of
the map $\Gamma \surj \pi_1(\widehat{M})$ factors through
$\pi_1(T^2/{\ZZ_3})=\{1\}$. So, the map $\Gamma \surj
\pi_1(\widehat{M})$ factors also through the quotient, $\ZZ^6
\surj \pi_1(\widehat{M})$. But $M$ contains three $2$--tori,
$T_1$, $T_2$ and $T_3$, which are the images of $\{
(x_1,x_2,0,0,0,0,0,0)\}$, $\{ (0,0,y_1,y_2,0,0,0,0)\}$ and $\{
(0,0,0,0,z_1,z_2,0,0)\}$, and $\pi_1(\widehat{M})$ is generated by
the images of $\pi_1(T_1)$, $\pi_1(T_2)$ and $\pi_1(T_3)$.
Clearly, $\ZZ_3$ acts in the standard way on each $T_i$. Therefore
$\pi_1(\widehat{M})$ is generated by $\pi_1(T_i/\ZZ_3)=\{1\}$,
which proves that $\pi_1(\widehat{M})=\{1\}$.
\end{proof}

\section{Non-formality of the quotient orbifold} \label{sec:non-formal}

Now we want to prove the non-formality of the orbifold
$\widehat{M}$ constructed in the previous section. By the results
of \cite{Ha, TO}, $M$ is non-formal since it is a nilmanifold
which is not a torus. We shall see that this property is inherited
by the quotient space $\widehat{M}=M/\ZZ_3$. For this, we study
the Massey products on $\widehat{M}$.

\begin{lemma}\label{lem:MasseyhatE}
 $\widehat{M}$ has a non-trivial Massey product  if and only if  $M$
 has a non-trivial Massey product with all cohomology classes
 $a_i\in H^*(M)$ being $\ZZ_3$-invariant cohomology classes.
\end{lemma}

\begin{proof}
We shall do the case of triple Massey products, since the general
case is similar. Suppose that $\la a_1,a_2,a_3\ra$,  $a_i \in
H^{p_i}(\widehat{M})$, $1 \leq i\leq 3$ is a non-trivial Massey
product on $\widehat{M}$. Let $a_i=[\alpha_i]$, where $\alpha_i\in
\Omega^*(\widehat{M})$. We pull-back the cohomology classes
$\alpha_i$ via $\varphi^*:\Omega^*(\widehat{M})\to \Omega^*(M)$ to
get a Massey product $\la
[\varphi^*\a_1],[\varphi^*\a_2],[\varphi^*\a_3]\ra$. Suppose that
this is trivial on $M$, then $\varphi^*\alpha_1\wedge
\varphi^*\alpha_2= d \xi$, $\varphi^*\alpha_2\wedge
\varphi^*\alpha_3=d \eta$, with $\xi,\eta\in \Omega^*(M)$, and
$\varphi^*\alpha_1 \wedge \eta+(-1)^{p_1+1} \xi \wedge
\varphi^*\alpha_3= df$. Then $\tilde\eta=
(\eta+\rho^*\eta+(\rho^*)^2\eta)/3$, $\tilde\xi=
(\xi+\rho^*\xi+(\rho^*)^2\xi)/3$ and $\tilde{f}= (f+
\rho^*\eta+(\rho^*)^2\eta)/3$ are $\ZZ_3$-invariant and
$\varphi^*\alpha_1 \wedge \tilde\eta+(-1)^{p_1+1} \tilde\xi \wedge
\varphi^*\alpha_3= d\tilde{f}$.

Conversely, suppose that $\la a_1,a_2,a_3\ra$, $a_i \in
H^{p_i}(M)^{\ZZ_3}$, $1 \leq i\leq 3$, is a non-trivial Massey
product on $M$. Then we can represent $a_i=[\alpha_i]$ by
$\ZZ_3$-invariant differential forms $\alpha_i\in
\Omega^{p_i}(M)$. Let $\hat\alpha_i$ be the induced form on
$\widehat{M}$. Then $\la [\hat\a_1],[\hat\a_2],[\hat\a_3]\ra$ is a
non-trivial Massey product on $\widehat{M}$. For if it were
trivial then pulling-back by $\varphi$, we would get $0\in \la
\varphi^*[\hat\a_1], \varphi^*[\hat\a_2],
\varphi^*[\hat\a_3]\ra=\la a_1,a_2,a_3\ra$.
\end{proof}

\begin{remark} \label{rem:sformal}
As $M$ is a nilmanifold which is not a torus, by \cite[Lemma
2.6]{FM1}, it is not $1$--formal. On the other hand, $\widehat{M}$
is simply connected by Proposition \ref{prop:thm_simply
connected0}, and hence it is $2$--formal. By the results of
\cite{FM1}, since $\widehat{M}$ is of dimension $8$, the only
possibility that it be non-formal is not to be $3$--formal. This
means that we have to compute the minimal model up to degree $3$,
which is a lengthy task, given that $b_2(\widehat{M})=13$ is quite
large. Therefore it is more convenient to find a suitable
non-trivial Massey product.
\end{remark}

\medskip

In our case, all the triple and quintuple Massey products on
$\widehat{M}$ are trivial. For instance, for a Massey product of
the form $\la a_1,a_2,a_3\ra$, all $a_i$ should have even degree,
since
$H^1(\widehat{M})=H^3(\widehat{M})=H^5(\widehat{M})=H^7(\widehat{M})=0$.
Therefore the degree of the cohomology classes in $\la
a_1,a_2,a_3\ra$ is odd, hence they are zero.

Since the dimension of $\widehat{M}$ is $8$, there is no room for
sextuple Massey products or higher, since the degree of $\la
a_1,a_2,\ldots, a_s\ra$ is at least $s+2$, as $\deg a_i\geq 2$.
For $s=6$, a sextuple Massey product of cohomology classes of
degree $2$ would live in the top degree cohomology. For computing
an element of $\la a_1,\ldots, a_6\ra$, we have to choose
$\alpha_{i,j}$ in (\ref{eqn:gm}). But then adding a closed form
$\phi$ with $a_1 \cup [\phi]=\lambda [\widehat{M}]\in
H^8(\widehat{M})$ to $\alpha_{2,6}$ we can get another element of
$\la a_1,\ldots, a_6\ra$ which is the previous one plus $\lambda
[\widehat{M}]$. For suitable $\lambda$ the we get $0\in \la
a_1,\ldots, a_6\ra$.

The only possibility for checking the non-formality of
$\widehat{M}$ via Massey products is to get a non-trivial
quadruple Massey product.

>From now on, we will denote by the same symbol a $\ZZ_3$-invariant
form on $M$ and that induced on $\widehat{M}$.  Notice that the
$2$ forms $\gamma_1\wedge\gamma_2$, $\beta_1\wedge\beta_2$ and
$\alpha_1\wedge\gamma_1+\alpha_2\wedge\gamma_1+\alpha_2\wedge\gamma_2$
are $\ZZ_3$-invariant forms on $M$, hence they descend to the
quotient $\widehat{M}=M/\ZZ_3$. We have the following:

\begin{proposition}\label{prop:nonformal1}
 The quadruple Massey product
 $$
 \la [\gamma_1\wedge\gamma_2],[\beta_1\wedge\beta_2],
 [\beta_1\wedge\beta_2],
 [\alpha_1\wedge\gamma_1+\alpha_2\wedge\gamma_1+\alpha_2\wedge\gamma_2]\ra
 $$
 is non-trivial on $\widehat{M}$.
 Therefore, the space $\widehat{M}$ is non-formal.
\end{proposition}

\begin{proof}
First we see that
 \begin{eqnarray*}
(\gamma_1\wedge\gamma_2)\wedge (\beta_1\wedge\beta_2)&=&d\xi, \\
(\beta_1\wedge\beta_2)\wedge
(\alpha_1\wedge\gamma_1+\alpha_2\wedge\gamma_1+\alpha_2\wedge\gamma_2)&=& d\varsigma,
 \end{eqnarray*}
 where $\xi$ and $\varsigma$ are the differential
$3$--forms on $\widehat{M}$ given by
  \begin{eqnarray*}
  \xi &=& - \frac{1}{6} \left( \gamma_1\wedge (\beta_1\wedge\eta_2+
   \beta_2\wedge\eta_2+\beta_2\wedge\eta_1 ) +\gamma_2\wedge(
   \beta_1\wedge\eta_2 + \beta_1\wedge\eta_1+ \beta_2\wedge\eta_1)\right), \\
  \varsigma &=&
  \frac{1}{3}\left(- \alpha_1\wedge(\eta_2\wedge\beta_1+
  \eta_1\wedge\beta_1+\eta_1\wedge\beta_2)+
  \alpha_2\wedge(\eta_2\wedge\beta_2-\eta_1\wedge\beta_1)\right).
 \end{eqnarray*}
Therefore, the triple Massey products $\la
[\gamma_1\wedge\gamma_2],[\beta_1\wedge\beta_2],
[\beta_1\wedge\beta_2]\ra$ and $\la [\beta_1\wedge\beta_2],
[\beta_1\wedge\beta_2],
[\alpha_1\wedge\gamma_1+\alpha_2\wedge\gamma_1+\alpha_2\wedge\gamma_2]\ra$
are defined, and they are trivial because all the (triple) Massey
products on $\widehat{M}$ are trivial. (Notice that the forms
$\xi$ and $\varsigma$ are $\ZZ_3$-invariant on $M$ and so descend
to $\widehat{M}$.) Therefore,  the quadruple Massey product $\la
[\gamma_1\wedge\gamma_2],[\beta_1\wedge\beta_2],
[\beta_1\wedge\beta_2],
[\alpha_1\wedge\gamma_1+\alpha_2\wedge\gamma_1+\alpha_2\wedge\gamma_2]\ra$
is defined on $\widehat{M}$. Moreover, it is trivial on
$\widehat{M}$ if and only if there are differential forms $f_i \in
\Omega^3(\widehat{M})$, $1\leq i\leq 3$, and $g_j \in
\Omega^4(\widehat{M})$, $1\leq j\leq 2$, such that
  \begin{eqnarray*}
  && (\gamma_1\wedge\gamma_2)\wedge (\beta_1\wedge\beta_2)=d(\xi+f_1), \\
  && (\beta_1\wedge\beta_2)\wedge (\beta_1\wedge\beta_2)=d f_2, \\
  && (\beta_1\wedge\beta_2)\wedge
   (\alpha_1\wedge\gamma_1+\alpha_2\wedge\gamma_1+\alpha_2\wedge\gamma_2)=
   d (\varsigma+f_3),  \\
  && (\gamma_1\wedge\gamma_2) \wedge f_2-(\xi+f_1)\wedge
   (\beta_1\wedge\beta_2)=d g_1,  \\
  && (\beta_1\wedge\beta_2)\wedge (\varsigma +f_3)-
   f_2\wedge (\alpha_1\wedge\gamma_1+
   \alpha_2\wedge\gamma_1+\alpha_2\wedge\gamma_2)=d g_2,
  \end{eqnarray*}
and the $6$--form given by
 $$
 \Psi=-(\gamma_1\wedge\gamma_2)\wedge g_2 - g_1\wedge
  (\alpha_1\wedge\gamma_1+\alpha_2\wedge\gamma_1+\alpha_2\wedge\gamma_2)
  +(\xi+f_1)\wedge(\varsigma +f_3)
 $$
defines the zero class in $H^6(\widehat{M})$. Clearly $f_1$, $f_2$
and $f_3$ are closed $3$--forms. Since $H^3(\widehat{M})=0$, we
can write $f_1=df'_1$, $f_2=df'_2$ and $f_3=df'_3$  for some
differential $2$--forms $f'_1$, $f'_2$  and $f'_3\in
\Omega^2(\widehat{M})$. Now, multiplying $[\Psi]$ by the
cohomology class $[\sigma]\in H^2(\widehat{M})$, where
$\sigma=2\alpha_1\wedge \gamma_2 -\alpha_2\wedge\gamma_1
+\alpha_1\wedge\gamma_1+\alpha_2\wedge\gamma_2$ we get
 $$
   \sigma \wedge \Psi
 =-\frac{1}{3}(\alpha_1\wedge\alpha_2\wedge\beta_1\wedge\beta_2
  \wedge\gamma_1\wedge\gamma_2\wedge\eta_1\wedge\eta_2)
   + d( \sigma\wedge\xi\wedge f_3' + \sigma \wedge \varsigma \wedge
 f_1'+\sigma \wedge  f_1' \wedge df_3').
 $$
Hence, $[2\alpha_1\wedge \gamma_2 -\alpha_2\wedge\gamma_1
+\alpha_1\wedge\gamma_1+\alpha_2\wedge\gamma_2]\cup [\Psi] \neq
0$, which implies that $[\Psi]$ is non-zero in $H^6(\widehat{M})$.
This proves that the Massey product $\la
[\gamma_1\wedge\gamma_2],[\beta_1\wedge\beta_2],
[\beta_1\wedge\beta_2],
[\alpha_1\wedge\gamma_1+\alpha_2\wedge\gamma_1+\alpha_2\wedge\gamma_2]\ra$
is non-trivial, and so $\widehat{M}$ is non-formal.
\end{proof}

\medskip

Let us see that $\widehat{M}$ is non-formal by proving that it has
a non-zero G-Massey product.

\begin{proposition} \label{prop:nonformal2}
Consider the following closed $2$--forms on $\widehat{M}$
 $$
\vartheta = \beta_1\wedge \beta_2, \quad
 \tau_1 = 2\alpha_1\wedge \gamma_2 -\alpha_2\wedge\gamma_1
  +\alpha_1\wedge\gamma_1+\alpha_2\wedge\gamma_2, \quad
 \tau_2 = \gamma_1\wedge\gamma_2,\quad
 \tau_3 = \alpha_1\wedge\gamma_1+\alpha_2\wedge\gamma_1+\alpha_2\wedge\gamma_2 .
 $$
Then the G-Massey product $\la [\vartheta];
[\tau_1],[\tau_2],[\tau_3] \ra$ is non-trivial on $\widehat{M}$.
\end{proposition}

\begin{proof}
A direct calculation shows that
 $$
\vartheta \wedge \tau_1 = d \kappa,
\quad
\vartheta \wedge \tau_2 = d \xi,
\quad
\vartheta \wedge \tau_3 = d \varsigma,
 $$
where $\xi$ and $\varsigma$ are the $3$--forms given in the proof
of Proposition
\ref{prop:nonformal1}, and $\kappa$ is the $3$--form
 \begin{eqnarray*}
 \kappa &=&\frac{1}{3}(
 \alpha_1\wedge \beta_1\wedge \eta_1-\alpha_1\wedge \beta_1\wedge \eta_2
 -\alpha_1\wedge \beta_2\wedge \eta_1-2\alpha_1\wedge \beta_2\wedge \eta_2  \\
 & &-\alpha_2\wedge \beta_1\wedge \eta_1-2\alpha_2\wedge \beta_1\wedge \eta_2
-2\alpha_2\wedge \beta_2\wedge \eta_1-\alpha_2\wedge \beta_2\wedge \eta_2).
\end{eqnarray*}

We know that the forms $\xi$ and $\varsigma$ are $\ZZ_3$-invariant
on $M$, and one can check that the form $\kappa$ is also. The
G-Massey product $\la [\vartheta]; [\tau_1],[\tau_2],[\tau_3] \ra$
is defined. As $H^3(\widehat{M})=0$, we have that $W=0$ in Lemma
\ref{lem:GMassey}, so $\la [\vartheta]; [\tau_1],[\tau_2],[\tau_3]
\ra$ consists of one element. This is
 $$
 \left[- \frac{4}{3} (\alpha_1\wedge\alpha_2\wedge\beta_1\wedge\beta_2
  \wedge\gamma_1\wedge\gamma_2\wedge\eta_1\wedge\eta_2)\right]
 \in H^8(\widehat{M}),
 $$
which is non-zero. So  $\la [\vartheta];
[\tau_1],[\tau_2],[\tau_3] \ra$ is non-trivial.
\end{proof}

\medskip

\section{Symplectic resolution of singularities} \label{sec:mainresult}

In this section we resolve symplectically the singularities of
$\widehat {M}$ to produce a {\em smooth\/} symplectic
$8$--manifold $\widetilde{M}$ which is simply connected and
non-formal. For this, we need the two following results:

\begin{proposition}\label{prop:symplectic}
The $2$--form $\omega$ on $M$ defined by
 $$
  \omega=\alpha_1\wedge\alpha_2+\eta_2\wedge\beta_1-
  \eta_1\wedge\beta_2+\gamma_1\wedge \gamma_2
 $$
is a $\ZZ_{3}$-invariant symplectic form on $M$. Therefore it
induces $\widehat\omega\in \Omega^2(\widehat{M})$, such that
$(\widehat{M},\widehat\omega)$ is a symplectic orbifold.
\end{proposition}

\begin{proof}
Clearly $\omega^4\not=0$. Using (\ref{eqn:v2}) we have that
$\rho^*(\omega)=(-\alpha_1-\alpha_2)\wedge\alpha_1+
\eta_1\wedge(-\beta_1-\beta_2)+(\eta_1+\eta_2)\wedge\beta_1+
(-\gamma_1-\gamma_2)\wedge\gamma_1 =\omega$, so $\omega$ is
$\ZZ_{3}$-invariant. Finally,
 $$
 d\omega= d\eta_2\wedge\beta_1- d\eta_1\wedge\beta_2
 =(\beta_2\wedge\gamma_1-\beta_2\wedge\gamma_2)\wedge\beta_1
 -(-\beta_1\wedge\gamma_1+\beta_1\wedge\gamma_2)\wedge\beta_2=
 0.
 $$
\end{proof}

\begin{lemma} \label{lem:coordinates}
Let $p \in M$ be a fixed point of the $\ZZ_3$-action. Then there
exists a system of complex coordinates $(w_1,w_2,w_3,w_4)$ around
$p$ such that the symplectic form $\omega$ defined in Proposition
\ref{prop:symplectic} can be expressed as
  $$
   \omega= i(dw_1\wedge d\bar{w}_1+dw_2\wedge d\bar{w}_2
   +dw_3\wedge d\bar{w}_3+  dw_4\wedge d\bar{w}_4 ).
  $$
Moreover, with respect to these coordinates, the $\ZZ_3$-action
$\rho$ on $M$ is given as
  $$
   \rho(w_1,w_2,w_3,w_4)=(\zeta^2 w_1,\zeta^2 w_2,\zeta w_3,\zeta^2 w_4),
  $$
where $\zeta=e^{\frac{2\pi i}{3}}$.
\end{lemma}

\begin{proof}
Let $p \in M$ be a fixed point of the $\ZZ_3$-action. Let $g\in G$
be a group element taking $p$ to the point $p_0=(0,\ldots, 0)\in
M$. Writing $m_g=m(g,\cdot)$, we have that $\rho \circ m_g=
m_{\rho(g)}\circ \rho$, so that $m_{\rho(g)}(p)=
\rho(m_g(p))=\rho(p_0)=p_0=m_g(p)$, therefore $\rho(g)$ coincides
with $g$ modulo $\Gamma$, and hence $\rho \circ m_g= m_g\circ
\rho$ on $M$. So we may suppose that the fixed point is $p=p_0$.

The coordinates for $G$ yield coordinates
$(x_1,x_2,y_1,y_2,z_1,z_2,v_1,v_2)$ for $M$ in a ball $B$ around
$p_0$ in which $p_0$ is mapped to the origin. The symplectic form
$\omega$ at the point $p_0$ is
 $$
  \omega_0=dx_1\wedge dx_2+d v_2\wedge dy_1
  -dv_1\wedge dy_2+dz_1\wedge dz_2.
 $$

Take now $\ZZ_3$-equivariant Darboux coordinates $\Phi \colon
(B,\omega) \longrightarrow (B_{\CC^4}(0,\epsilon),\omega_0)$, for
some $\epsilon>0$. This means that $\Phi\circ
d\rho_{p_0}=\rho\circ \Phi$ and $\Phi^*\omega_0=\omega$. The proof
of the existence of usual Darboux coordinates in \cite[pp.\
91--93]{McDuff-Salamon} carry over to this case, only being
careful that all the objects constructed should be
$\ZZ_3$-equivariant.

In terms of the complex coordinates $x=x_1+ix_2$, $y=y_1+iy_2$,
$z=z_1+iz_2$, $v=v_1+iv_2$ of $\CC^4$, the form $\omega$ is
written as
 $$
  \omega=i(dx\wedge d\bar{x}+dv\wedge d\bar{y}+
  dy\wedge d\bar{v}+dz\wedge d\bar{z}).
 $$
Now, we define the functions $u=\frac{1}{\sqrt{2}}(v+y)$ and
$w=\frac{1}{\sqrt{2}}(v-y)$. Since $dv\wedge d\bar{y}+dy\wedge
d\bar{v}=du\wedge d\bar{u}- dw\wedge d\bar{w}$, the symplectic
form $\omega$ is expressed as
 $$
  \omega=i(dx\wedge d\bar{x}+du\wedge d\bar{u}
  -dw\wedge d\bar{w}+dz\wedge d\bar{z}).
 $$
Consider new complex functions $x'=x'_1+ix'_2$, $u'=u'_1+iu'_2$,
$w'=w'_1+iw'_2$ and $z'=z'_1+iz'_2$, where
 $$
 \left\{
 \begin{array}{ll}
  x'_1=x_1-\frac{1}{2}x_2,   \quad & x'_2=-\frac{\sqrt{3}}{2}x_2, \\
  u'_1=u_1-\frac{1}{2}u_2,   \quad & u'_2=-\frac{\sqrt{3}}{2}u_2,
 \end{array}
 \right. \qquad
 \left\{
 \begin{array}{ll}
  w'_1=w_1-\frac{1}{2}w_2, \quad  & w'_2=-\frac{\sqrt{3}}{2}w_2, \\
  z'_1=z_1-\frac{1}{2}z_2, \quad  & z'_2=-\frac{\sqrt{3}}{2}z_2.
 \end{array}
 \right.
 $$
So,
 $$
  \rho (x',u',w',z')= (\zeta x', \zeta u',\zeta w',\zeta z'),
 $$
with $\zeta=e^{\frac{2\pi i}{3}}$ (by using that $\rho$
corresponds to the matrix in (\ref{eqn:matrix}) for the
coordinates $x,u,w,z$). Since $dx'\wedge
d\bar{x'}=\frac{-\sqrt{3}}{2} dx\wedge d\bar{x}$, we have
 $$
  \omega=-\frac{2i}{\sqrt{3}}(dx'\wedge d\bar{x'}+du'\wedge d\bar{u'}
  -dw'\wedge d\bar{w'}+dz'\wedge d\bar{z'}).
 $$
Finally, the set of coordinates $(w_1,w_2,w_3,w_4)=
\sqrt{\frac{2}{\sqrt{3}}} (\bar{x'},\bar{u'},w',\bar{z'})$ gives
the desired result.
\end{proof}

\medskip
Next, we see how it is possible to desingularize the space
$\widehat{M}$. We use the following result which is \cite[Lemma
2.2]{FM4}.

\begin{lemma}\label{lem:glue}
 Let $(B,\omega_0)$ be the standard K{\"a}hler ball in $\CC^n$, $n>1$, and let $\Pi$
 be a finite group acting linearly (by complex isometries)
 on $B$ whose only fixed point is the origin.
 Let $\phi:(\widetilde{B},\omega_1) \to (B/\Pi,\omega_0)$ be a K{\"a}hler
 resolution of the singularity of the quotient. Then there is a
 K{\"a}hler form on $\widetilde{B}$ such that it coincides with
 $\omega_0$ near the boundary, and with a positive multiple of
 $\omega_1$ near the exceptional divisor $E=\phi^{-1}(0)$. \hfill
 $\Box$
\end{lemma}

\begin{theorem}\label{thm:prop-desingularization}
 There is a smooth compact
 symplectic manifold $(\widetilde{M},\widetilde{\omega})$
 which is isomorphic to
 $(\widehat{M},\widehat{\omega})$ outside the singular points.
\end{theorem}

\begin{proof}
 Let $p$ be a fixed point of the $\ZZ_3$-action.
 By Lemma \ref{lem:coordinates} we have a K\"ahler model for a
 neighbourhood $B$ of $p$, where the action is of the form
 $(w_1,w_2,w_3,w_4) \mapsto
 (\zeta^2 w_1,\zeta^2 w_2,\zeta w_3,\zeta^2 w_4)$.
 We may resolve the singularity of $B/\ZZ_3$ with a K\"ahler
 model. We do the resolution of singularities via iterated
 blow-ups as it is a standard procedure for algebraic manifolds.

Blow up $B$ at $p$ to get $\widetilde{B}$. This replaces the point $p$
by a complex projective space $F=\PP^3$ in which $\ZZ_3$ acts as
  $$
   [w_1,w_2,w_3,w_4] \mapsto
   [\zeta^2 w_1,\zeta^2 w_2,\zeta w_3,\zeta^2 w_4]=
   [w_1,w_2,\zeta^2 w_3,w_4].
  $$
Therefore there are two components of the fix-point locus of the
$\ZZ_3$-action on $\widetilde{B}$, namely the point $q=[0,0,1,0]$
and the complex projective plane $H=\{[w_1,w_2,0,w_4]\} \subset
F=\PP^3$. Next blow up $\widetilde{B}$ at $q$ and at $H$ to get
$\widetilde{\widetilde{B}}$. The point $q$ is substituted by a
projective space $H_1=\PP^3$. The normal bundle of $H\subset
\widetilde{B}$ is the sum of the normal bundle of $H\subset F$,
which is $\cO_{\PP^2}(1)$, and the restriction of the normal
bundle of $F\subset\widetilde{B}$ to $H$, which is
$\cO_{\PP^3}(-1)|_{\PP^2}= \cO_{\PP^2}(-1)$. Therefore the second
blow-up replaces the plane $H$ by the $\PP^1$-bundle over $\PP^2$
defined as $H_2 = \PP( \cO_{\PP^2}(-1)\oplus\cO_{\PP^2}(1))$. The
strict transform of $F\subset\widetilde{B}$ under the second
blow-up is the blow up $\widetilde{F}$ of $F=\PP^3$ at $q$, which
is a $\PP^1$-bundle over $\PP^2$, actually $\widetilde{F} = \PP(
\cO_{\PP^2}\oplus\cO_{\PP^2}(1))$.

The fix-point locus of the $\ZZ_3$-action on
$\widetilde{\widetilde{B}}$ consists of the two disjoint divisors
$H_1$ and $H_2$. Therefore the quotient
$\widetilde{\widetilde{B}}/\ZZ^3$ is a smooth K\"ahler manifold
\cite[page 82]{BPV}. This provides a symplectic resolution of the
singularity $B/\ZZ_3$. To glue this K\"ahler model to the
symplectic form in the complement of the singular point using
Lemma \ref{lem:glue}. We do this at every fixed point to get a
smooth symplectic resolution of $\widehat{M}$.
\end{proof}

%
%
%

\begin{theorem}\label{thm::nonformal3}
 The manifold $\widetilde{M}$ is non-formal.
\end{theorem}

\begin{proof}
All the forms of the proof of either Proposition
\ref{prop:nonformal1} or Proposition \ref{prop:nonformal2} can be
defined on the resolution $\widetilde{M}$ as follows: take a
$\ZZ_3$-equivariant map $\psi:M\to M$ which is the identity
outside small balls around the fixed points, and contracts smaller
balls onto the fixed points. Substitute the forms $\vartheta$,
$\tau_i$, $\kappa$, $\xi$, \ldots\ by $\psi^*\vartheta$,
$\psi^*\tau_i$, $\psi^*\kappa$, $\psi^*\xi$, \ldots\ Then the
corresponding elements in the quadruple Massey product or the
G-Massey product are non-zero, but these forms are zero in a
neighbourhood of the fixed points. Therefore they define forms on
$\widetilde{M}$, by extending them by zero.
\end{proof}

\begin{theorem}\label{thm::simply connected2}
 The manifold $\widetilde{M}$ is simply connected.
\end{theorem}

\begin{proof}
We have already seen in Proposition \ref{prop:thm_simply
connected0} that $\widehat{M}$ is simply connected. The resolution
$\widetilde{M}\to \widehat{M}$ consists of substituting, for each
singular point $p$, a neighbourhood $B/\ZZ_3$ of it by the
non-singular model $\widetilde{\widetilde B}/\ZZ_3$. The fiber
over the origin of $\widetilde{\widetilde B}/\ZZ_3\to B/\ZZ_3$ is
simply connected: it consists of the union of the three divisors
$H_1=\PP^3$, $H_2=\PP( \cO_{\PP^2}(-1)\oplus\cO_{\PP^2}(1))$ and
$\widetilde F/\ZZ_3=\PP( \cO_{\PP^2}\oplus\cO_{\PP^2}(3))$, all of
them are simply connected spaces, and their intersection pattern
forms no cycles. A Seifert-Van Kampen argument proves that
$\widetilde{M}$ is simply connected.
\end{proof}

\begin{remark}
 The second Betti number of $\widetilde{M}$ increases in $3$ by
 a desingularisation of a fixed point. As there are $81$ fixed
 points, we have $b_2(\widetilde{M})=b_2(\widehat{M})+81\cdot
 3=256$. This makes very difficult to write down the minimal model
 of $\widetilde{M}$ up to degree $3$ to check non-formality.
\end{remark}

\begin{remark} \label{rem:hardlefschetz}
The symplectic orbifold $(\widehat{M},\widehat\omega)$ is not hard
Lefschetz: consider the $\ZZ_3$-invariant forms $\beta_1\wedge
\beta_2$ and $\alpha_1\wedge\alpha_2\wedge \xi$ on $M$, where
$\xi$ is the $3$--form given in the proof of Proposition
\ref{prop:nonformal1}. As they are $\ZZ_3$-invariant forms, they
descend to $\widehat{M}$. But then
 $$
 \omega^2 \wedge (\beta_1\wedge
 \beta_2)= \alpha_1\wedge\alpha_2\wedge\beta_1\wedge \beta_2\wedge
 \gamma_1\wedge\gamma_2= d(\alpha_1\wedge\alpha_2\wedge \xi),
 $$
which means that the map $[\omega]^2 \colon H^2(\widehat{M}) \to
H^6(\widehat{M})$ is not a monomorphism. These forms can be
extended to $\widetilde{M}$ via the process carried out in the
proof of Theorem \ref{thm::nonformal3}. Therefore, the map
$[\omega]^2 \colon H^2(\widetilde{M}) \to H^6(\widetilde{M})$ is
not injective.

This raises the question of the existence of a non-formal simply
connected compact symplectic $8$--manifold satisfying the hard
Lefschetz property (Cavalcanti \cite{Cav2} has given examples for
dimensions $\geq 10$).
\end{remark}

{\small

\vspace{0.15cm}

\noindent{\sf M. Fern\'andez:} Departamento de Matem\'aticas,
Facultad de Ciencia y Tecnolog\'{\i}a, Universidad del Pa\'{\i}s
Vasco, Apartado 644, 48080 Bilbao, Spain.

{\sl E-mail:} marisa.fernandez@ehu.es

\vspace{0.15cm}

\noindent{\sf V. Mu\~noz:} Departamento de Matem\'aticas, Consejo
Superior de Investigaciones Cient{\'\i}ficas, C/ Serrano 113bis, 28006
Madrid, Spain.

Facultad de Matem{\'a}ticas, Universidad Complutense de Madrid, Plaza
de Ciencias 3, 28040 Madrid, Spain.

{\sl E-mail:} vicente.munoz@imaff.cfmac.csic.es}

\end{document}